\documentclass[12pt, reqno]{amsart}
\usepackage[margin=2cm]{geometry}
\usepackage[utf8]{inputenc}
\usepackage{comment}
\usepackage{hyperref}
\usepackage{amsfonts}
\usepackage{amsmath}
\usepackage{amssymb}
\usepackage{amscd}
\usepackage{graphicx}
\usepackage{latexsym}
\usepackage{color}
\usepackage{epsfig}
\usepackage{amsopn}
\usepackage{xypic}
\usepackage{mathrsfs}
\usepackage{array}
\usepackage{xspace}
\usepackage{comment}
\usepackage[hypcap]{caption}
\usepackage{caption}
\usepackage{subcaption}
\usepackage[colorinlistoftodos]{todonotes}
\usepackage{verbatim}
\usepackage{enumitem}
\usepackage{amsaddr}

\newtheorem{thm}{Theorem}[section]
\newtheorem{lem}[thm]{Lemma}

\newtheorem{prop}[thm]{Proposition}

\theoremstyle{definition}
\newtheorem{definition}[thm]{Definition}

\newtheorem{exam}[thm]{Example}

\theoremstyle{remark}
\newtheorem{rmk}[thm]{Remark}

\newcommand{\C}{\mathbb{C}}
\newcommand{\Z}{\mathbb{Z}}

\newcommand{\RR}{\mathbb{R}}

\newcommand{\To}{\longrightarrow}

\newcommand{\cDD}{\mathcal{D}}
\newcommand{\cP}{\mathcal{P}}
\newcommand{\cQ}{\mathcal{Q}}

\newcommand{\cA}{\mathcal{A}}

\newcommand{\Aut}{\text{Aut}}
\newcommand{\Hom}{\mathrm{Hom}}
\newcommand{\Stab}{\mathrm{Stab}}

\newcommand{\GStab}{\mathrm{GStab}}
\newcommand{\HN}{\mathrm{HN}}

\title{A local compactification of the Bridgeland stability manifold}
\author{Barbara Bolognese}
\thanks{This project was partially supported by the ERC Advanced Grant {\it Stability conditions, Donaldson-Thomas invariants and cluster varieties}, P.I. Tom Bridgeland.}
\address{Department of Mathematics, Universit\`a degli studi di Roma Tre}
\email{barbara.bolognese@uniroma3.it}

\begin{document}

\maketitle

\begin{abstract}
Bridgeland stability manifolds of Calabi-Yau categories are of noticeable interest both in mathematics and in physics. By looking at some of the known example, a pattern clearly emerges and gives a fairly precise description of how they look like. In particular, they all seem to have missing loci, which tend to correspond to degenerate stability conditions vanishing on spherical objects. Describing such missing strata is also interesting from a mirror-symmetric perspective, as they conjecturally parametrize interesting types of degenerations of complex structures. All the naive attempts at constructing modular partial compactifications show how elusive and subtle the problem in fact is: ideally, the missing strata would correspond to stability manifolds of quotient triangulated categories, but establishing such correspondence on geometric level and viewing stability conditions on quotients of the original triangulated category as suitable degenerations of stability conditions is not straightforward. In this paper, we will present a method to construct such partial compactifications if some additional hypoteses are satisfied, by realizing our space of interest as a suitable metric completion of the stability manifold. 
\end{abstract}

\section{Introduction}

One of the most remarkable features of the theory of Bridgeland stability is the statement, originally proved by Bridgeland in \cite{Bri}, that the set of Bridgeland stability conditions on a given triangulated category $\cDD$, if non empty, has the structure of a complex manifold, hence it is usually referred to as the {\it stability manifold} of $\cDD$. The complex structure is induced by proving that the forgetful map, sending a stability condition consisting of a pair of a central charge and a slicing to its central charge, is a local homeomorphism onto an open subset of a complex vector space. Stability manifolds are central in the most recent development of various branches of mathematics, such as mirror symmetry, geometric representation theory and algebraic geometry, and have been completely described in various classes of examples. Variation of stability shows that there is a wall-and-chamber decomposition within the stability manifold: proper semi-stability tends to occur on walls, i.e. real codimension one submanifolds, and stable objects tend to stay stable when the stability condition varies within a chamber. This phenomenon often induces an interesting behavior in the corresponding moduli spaces of (semi)-stable objects: one observes that, in most cases, crossing a wall leads to a birational transformation between the moduli spaces of stable objects in the two adjacent chambers. This has led to a detailed description of the birational geometry of classical moduli spaces of stable sheaves in several cases, as we are about to explain.

In algebraic geometry, i.e. when one considers the stability manifold of the derived category of the derived category of a (smooth) projective variety, the theory is able to give a complete description of the stability manifold of curves (\cite{Bri} for elliptic curves,  \cite{Ok} and \cite{BSW} for $\mathbb{P}^1$ and \cite{Mac} for curves of genus at least two), K3 surfaces (\cite{BriK3}), \cite{BB}), the local $\mathbb{P}^2$ (\cite{BMplane}) and there are various partial descriptions of the wall-and-chamber decomposition on smooth surfaces in general (\cite{Mac}). In geometric representation theory, one usually considers two different triangulated categories associated to a quiver: the derived category of its path algebra, and the derived category of its Ginzburg algebra. Examples in this setting include the $A_2$ quiver (\cite{BSQ}), the Kronecker quiver (\cite{McA}), the acyclic triangular quiver (\cite{DK}) and a complete description for quivers coming from the triangulation of surfaces has been given in \cite{BS}, relating stability manifolds to moduli spaces of quadratic differential of Riemann surfaces. 

\

In virtue of this, it becomes interesting to study degenerations of stability conditions. The mirror-symmetric conjectural interpretation of the stability manifolds (see e.g. \cite{Brispace}, \cite{HW}) suggests that two types of degenerations are of particular interest.
The first type, called {\it large volume limit}, arises in the algebro-geometric setting when considering stability conditions which corresponds to a ``large" polarization. Stability conditions in the large volume limit are actual stability conditions and identify a special chamber, called the {\it Gieseker chamber}, in which stable complexes are in fact concentrated in one degree. Via the birational maps arising from the wall-and-chamber decomposition described in the previous paragraph, this gives some remarkable information about the birational geometry of moduli spaces of sheaves (see e.g. \cite{Yos}, \cite{BM} \cite{BM2}, \cite{ABCH}, \cite{Nu}, \cite{BHLRSWZ}).
The second type, coming from a more mirror-symmetric approach, arises when considering stability condition acquiring a massless object in the limit. Those should correspond to maximally unipotent degenerations of K\"ahler structures, and in the space of central charges are usually parametrized by subvarieties which are not in the image of the forgetful map and appear as ``missing" loci in the space of central charges.  Degenerate stability conditions in this sense are not stability conditions, as the requirement that the mass of every object be positive fails to be satisfied, although they still have a viable categorical interpretation as we are about to explain. The purpose of this paper is that of describing a partial compactification of the stability manifolds having degenerate stability conditions of this type appearing as boundary strata. One might look at the problem on two different levels: the topological level, and the categorical level.

On a topological level, it appears natural to give a description of degenerate stability conditions as limits of Cauchy sequences and consider the metric completion of the stability manifold. Indeed this is what seems to happen in the space of central charges: degenerate central charges in this sense appear in the closure of the image of the forgetful map, which coincides with its metric completion. To this purpose, the natural metric introduced by Bridgeland on the set of stability conditions does not seem to be an effective choice, since one can show that stability conditions acquiring massless objects appear infinitely distant to any stability condition. Indeed, one can show that the stability manifold with this metric is already complete (\cite{Woo}). Hence we introduce a new metric, constructed as the geodesic metric associated to the pullback of any inner product from the space of central charges via the forgetful map. We consider the metric completion $\widehat{\Stab}(\cDD)$ of the stability manifold with respect to this metric. Such a metric completion seems to produce a meaningful geometric object which serves the purpose on a topological level, even though it does not offer any modular interpretation of the objects appearing in the added strata.

On a categorical level, we would like to identify degenerate stability conditions with proper stability conditions on a quotient triangulated category. We define {\it the set of generalized stability conditions} as the set of pairs 
\[ \text{GStab}(\cDD)= \{ (\mathcal{K}, \sigma) \ | \ \mathcal{K} \subset \mathcal{D} \text{ thick subcategory }, \ \sigma \in \Stab(\mathcal{D}/\mathcal{K})\}  \]

and, in view of the natural injection $\Stab(\mathcal{D}) \hookrightarrow \text{GStab}(\cDD)$ given by sending a stability condition $\sigma$ to the pair $(0, \sigma) \in \text{GStab}(\cDD)$, declare it to be the natural local compactification of the stability manifold $\Stab(\cDD)$. There are two main issues that one can easily find: first of all, the set $\text{GStab}(\cDD)$ does not seem to have a natural topology, as there seems not to be an easy way of viewing the loci which correspond to pairs where $\mathcal{K}\neq 0$ as boundary strata of the stability manifold. Secondly, the set of generalized stability conditions seems unnaturally bigger than the set of stability conditions, as it would not be a reasonable expectation to have any thick subcategory in $\cDD$ appear in relation to a degenerate stability condition. 

Our main result is to show that these two construction are relatable to each other. We prove the following:

\begin{thm}\label{Main}
Let $\cDD$ be a triangulated category, and assume that the natural forgetful morphism 
\[ \pi : \Stab(\cDD) \to \Hom _{\Z}(K({\cDD}), \mathbb{C}) , \ (Z, \mathcal{P}) \mapsto Z \]
is a covering map onto the open dense submanifold $X = \Hom_{\Z}(K(\cDD), \mathbb{C}) \setminus \Delta\subset \Hom _{\Z}(K(\cDD), \mathbb{C})$, where $\Delta$ is a locally finite collection of smooth submanifolds. Then
there is an injective morphism of sets \[ j: \widehat{\Stab(\cDD) } \longrightarrow \GStab(\cDD) . \] Its image $\overline{\Stab(\cDD)}=j(\widehat{\Stab(\cDD)})\subset \GStab(\cDD)$, called the set of {\it degenerate} stability conditions, thus inherits a natural topology. 
\end{thm}

It might be interesting to point out that our construction seems a generalization of that described in \cite{All}, where the author describes a particular metric completion space and proves that it is a CAT(0) space. CAT(0) spaces are known to be contractible (see e.g. \cite{Cap}), hence it might be interesting to explore the question whether the same results hold for the stability manifold with the construction we propose in our paper. A few consideration on CAT(0) properties of stability manifolds have been made in \cite{Kik}.

The paper will be structured as follows. Section (\ref{Preliminaries}) will recall a few generalities on Bridgeland stability conditions. In section (\ref{lms}), we will introduce the notion of a locally modelled metric space and study its properties. Section (\ref{pbm}) will realize certain stability manifolds with a suitable metric as locally modelled spaces. In section (\ref{tls}) and (\ref{lddsc}) we will construct the map $j$ and prove Theorem \ref{Main}. 

\subsection*{Acknowledgements} The author wishes to heartily thank Tom Bridgeland: this paper is the result of extensive and frequent conversations with him, and has much benefited from his ideas, his knowledge and his suggestions. The author would also like to thank Domenico Fiorenza for pointing out a few mistakes in an early version of the manuscript.
\section{Preliminaries}\label{Preliminaries}

We will now recall a few preliminaries on Bridgeland stability conditions: for a more complete treatment, see \cite{Bri}. Let $\mathcal{D}$ be a triangulated category, and let $K(\cDD)$ be its Grothendieck group. A {\it Bridgeland stability condition} (BSC) on $\cDD$ is given by a pair $\sigma = (Z, \cP)$, where:
\begin{enumerate}
	\item the {\it central charge} $Z: K(\cDD) \To \C $ is a group homomorphism;
	\item the {\it slicing} $\cP = \{ \cP(\phi) \}_{\phi \in \RR}$ is a collection of full additive subcategories $\cP(\phi)$ satisfying:
	\begin{itemize}
		\item[(a)] For every $\phi \in \RR$, one has $\cP (\phi + 1) = \cP (\phi) [1]$;
		\item[(b)] For every $\phi _1 > \phi _2$, one has $\Hom( \cP(\phi_1), \cP(\phi _2))=0$;
		\item[(c)] For every object $E\in \cDD$ there exists a Harder-Narasimhan filtration, i.e. a (unique) collection of distinguished triangles (DTs):
		$$
		\hspace{-0.2in}
		\small{
		\xymatrix{
		0=E_0 \ar[rr] & & E_1 \ar[rr] \ar[dl] & &  \cdots \ar[rr]  \ar[dl] & & E_{n-1} \ar[rr]  \ar[dl] & & E_n=E  \ar[dl] \\
		& F_1 \ar[lu]^{+1}&    &F_2 \ar[lu]^{+1} &    &F_{n-1} \ar[lu]^{+1} &   &F_n \ar[lu]^{+1} & \\
		}}
		$$
		where $F_i\in \cP (\phi _i)$ for each $i=1,...,n$ and $\phi _1 > \phi _2 > \cdots > \phi _n$;
	\end{itemize}
	\item For each nonzero $E\in \cP (\phi) $ for $\phi \in \RR$, one can write $Z(E) = r(E) e^{i\pi \phi}$ with $r(E) \in \RR _{>0}$.
\end{enumerate}
For an object $E\in \cDD $, we will refer to the real numbers $\phi _1, ... , \phi _n$ obtained above as the {\it phases} of $E$. Moreover, we will set $\phi ^+ _{\cP } (E) = \phi _1$ (resp. $\phi ^- _{\cP } (E) = \phi _n$), which will be henceforth called the {\it maximal phase} (resp. {\it minimal phase}) of $E$. Similarly, if $\sigma$ is a stability condition, we will denote by $\phi ^+ _{\sigma } (E) = \phi _1$ (resp. $\phi ^- _{\sigma } (E) = \phi _n$) the maximal (resp. minimale) phase of $E$ with respect to the slicing associated to $\sigma$. Equivalently, one can give a BSC as a pair $(Z,\cA)$ where the central charge $Z$ is as above and $\cA$ is the heart of a bounded t-structure on $\cDD$, such that for each nonzero object $E\in \cA $, $Z(E)$ belongs to the upper half plane 
 \[\mathbb{H} = \{ z\in \C \ | \ \Im z \geq 0 \text{ and } \Im z = 0 \Rightarrow \Re Z < 0 \}. \]
 The equivalence of these two definitions is easily seen: if one sets
 \[ \cP (\geq \psi) = \{ E\in \cDD \ | \ \text{the phases of } E \text{ wrt } \cP \text{ lie in the interval } [\psi , +\infty ) \} \]
it is easy to show that the pair $( \cP (\geq \psi), \cP (< \psi + 1))$ gives a bounded t-structure on $\cDD $, whose heart is $\cA = \cP ([\psi, \psi + 1))$. Conversely, given a pair $(Z, \cA )$ as above, one can produce a slicing on $\cDD $ by combining the existence of two natural HN filtrations: the one on $\cA$ in terms of phases, and the one on $\cDD $ in terms of shifts of the heart $\cA$. We will use both definitions from now on.

For each object $E\in \cDD$ and each stability condition $\sigma$, we denote $m_{\sigma}(E)$ denotes the {\it mass} of $E$ with respect to $\sigma$, i.e. the sum of the amplitudes of its Harder-Narasimhan factors:
	\[ m_{\sigma}(E) = \sum _{i=1}^k |Z(\mathrm{HN}_{\sigma}^i(E))| . \]
Condition (3) above tells us that the mass of any object is strictly positive.

One can set :
\[ \text{Slice} (\cDD ) = \{ \cP \ | \ \cP \text{ is a slicing on }\cDD \} , \text{ and} \]
\[ \Stab (\cDD) = \{ \sigma = (Z, \cP ) \ | \ \sigma \text{ is a BSCs on }\cDD \} \subset \Hom _Z(K(\cDD ) , \C ) \times \text{Slice}(\cDD ). \]
A priori, $\Stab (\cDD)$ is a set without any additional structure. It is possible to construct a natural topology on $\Stab (\cDD)$ following the steps below:
\begin{enumerate}
	\item Construct a generalized norm on the complex vector space $\Hom _Z(K(\cDD ) , \C )$. Fix a stability condition $\sigma = (Z, \cP )$ in $\Stab (\cDD )$. For each $U \in \Hom _Z(K(\cDD ) , \C )$, one can set:
	\[ || U || _{\sigma } = \underset{E \ \sigma -\text{ss}}{\text{sup}} \frac{|U(E)|}{|Z(E)|} .\]
	It is easy to show that this defines a generalized norm, which is moreover indepentent on the stablity condition we fixed at the beginning.
	\item Construct a generalized metric on the set $\text{Slice} (\cDD )$. For each pair of slicings $\cP$ and $\cQ$ in $\text{Slice} (\cDD )$, set
	\[ d(\cP, \cQ) =  \underset{0\neq E \in \cDD}{\text{sup}} \{ |\phi ^-_{\cP}(E) - |\phi ^-_{\cQ}(E)|, |\phi ^+_{\cP}(E) - |\phi ^+_{\cQ}(E)| \}  .\]
	One can show that this indeed defines a generalized metric. 
	\item Give $\Stab (\cDD)$ a topology in terms of a basis of open balls. For each $\sigma = (Z, \cP ) \in \Stab (\cDD)$, define:
	\[ B_{\epsilon}(\sigma ) = \{ \tau = (W, \cQ ) \ | \ ||W-Z||_{\sigma } < \sin \epsilon \text{ and } d(\cP , \cQ)<\epsilon \} .\]
One can show that this is a basis of open subsets for a topology, see e.g. \cite[Lemma 6.2]{Bri}.
\end{enumerate}

\

Moreover, one can show that this topology is induced by a generalized metric on the whole $ \text{Stab} (\cDD )$: for each pair of BSCs $\sigma _1=(Z, \cP ), \sigma _2= (W, \cQ) \in  \text{Stab} (\cDD )$, one can set: 
	\[ d(\sigma _1, \sigma _2) =  \underset{0\neq E \in \cDD}{\text{sup}} \left\{ |\phi ^-_{\sigma _1}(E) - \phi ^-_{\sigma _2}(E)|, |\phi ^+_{\sigma _1}(E) - |\phi ^+_{\sigma _2}(E)| \} , \left|\log \frac{m_{\sigma _1}(E)}{m_{\sigma _2}(E)}\right| \right\}.\]
	
Consider a sequence of stability conditions $\sigma _n \in \Stab(\cDD)$, and suppose that there exists an object $E\in \cDD$ such that $\lim _{n\to \infty} m_{\sigma _n}(E) = 0$. Then for any stability condition $\tau$, one sees that
\[ \lim _{n\to \infty}d(\sigma _n, \tau) = +\infty . \]
We informally call a sequence of stability conditions acquiring a massless object at infinity a {\it degenerate} stability condition. Our observation shows that this metric does not realize degenerate stability conditions as missing strata in the stability manifold, but rather as strata at infinity. Indeed, one can show that the space of stability conditions  with this metric is complete, see \cite{Woo}. 	
	
From now on, we will only consider stability conditions $\sigma = (Z, \cP)$ satisfying the following additional assumptions:
\begin{enumerate}
	\item The slicing $\cP$ is {\it locally finite}, i.e. there exists $\eta \in \RR _{\geq 0}$ \ such that for each $\phi \in \RR $ the category $\cP ((\phi - \eta , \phi + \eta))$ is of finite length.
	\item Fix a finite dimensional lattice $\Lambda $ and a surjection $p: K(\cDD ) \to \Lambda$. Then the group homomorphism $Z$ must factor through $p$.
\end{enumerate}

With a slight abuse of notation, we will call $ \Stab (\cDD)$ the subspace of stability condition satisfying hypoteses (1)-(2) above. Then one has the following:
\begin{thm}[\cite{Bri}]\label{forgetful} The forgetful morphism
\[ \pi : \Stab (\cDD) \To \Hom _Z(K(\cDD ) , \C )\]
is a local homeomorphism, hence the space $\Stab (\cDD)$, if non empty, has the structure of a complex manifold.
\end{thm}

In order to see the degenerate stability conditions as missing points in some metric completion of the stability manifold, we will construct a new metric, which will have the important property that degenerate stability conditions will appear as limits of Cauchy sequences. The stability manifolds with this new metric will no longer be complete, and we will show that degenerate stability conditions will appear in the added strata of its metric completion. The new metric will be constructed as follows: consider the forgetful morphism $\pi$ appearing in Theorem (\ref{forgetful}). Fix an inner product $g$ on the real vector space underlying the complex vector space $\Hom _Z(K(\cDD ) , \C )$, look at this inner product as a translation invariant Riemannian metric on $\Hom _Z(K(\cDD ) , \C )$ and consider the {\it pullback} Riemannian metric $\tilde{g}= \pi^*g$ on $\Stab (\cDD)$. Since the map $\pi$ is a local homeomorphism, the pair $(\Stab (\cDD), \tilde{g})$ is a Riemannian manifold, hence it is a length space with respect to the induced Riemannian distance. 

\section{Locally modelled metric spaces}\label{lms}
In this section we introduce the main topological tools to construct the map $j$ appearing in Theorem \ref{Main}. It has a high technical content, hence it can be skipped or quickly read through by a reader who is mainly interested in its application to stability conditions. An interesting reference on length spaces and their properties can be found in \cite{BT}.
Let $(X, d)$ and $(Y, d')$ be metric spaces, and let $\pi: Y \to X$ be a local isometry with respect to the two topologies induced by the two metrics. We are interested in studying the metric completion of the source space $Y$ and in relating it to that of the target space $X$. With this purpose in mind, we define and study a special class of Cauchy sequences on $Y$ which behave particularly favorably with respect to the local homeomorphism $\pi$.

\begin{definition}\label{good}
A Cauchy sequence $\underline{y}=\{ y_n \} $ in $Y$ is called {\it $\pi$-local} if there exists $N>0$ and $U\subset Y$ open subset such that
\begin{enumerate}
\item $\{ y_n \} _{n>N} \subset U$
    \item  $\pi |_{U} : U \to \pi(U)$ is a homeomorphism.
\end{enumerate}
\end{definition}

For a $\pi$-local Cauchy sequence $\underline{ y }$, we denote the choice of such an open subset by $U_{\underline{y}}$. The notion of $\pi$-local Cauchy sequence turns out to be surprisingly convenient in this set-up. We define a new equivalence relation on $\pi$-local Cauchy sequences: 

\begin{definition}
Two $\pi$-local Cauchy sequences $\underline{y }$ and $\underline{y}'$ in $Y$ are called {\it very strongly equivalent} if one can choose $U_{\underline{y }}=U_{\underline{y} '}$.
\end{definition}

Notice that being very strongly equivalent is not an equivalence relation, since it is not transitive. Hence, we give the following:

\begin{definition}
Two $\pi$-local Cauchy sequences $\underline{y }$ and $\underline{y} '$ in $Y$ are called {\it strongly equivalent} if there is a chain $$\underline{y}= \underline{y^0} \sim \underline{y^1} \sim \underline{y^2} \sim ... \sim \underline{y^m}= \underline{y}'$$ with $\underline{y ^i}$ very strongly equivalent to $\underline{y ^{i+1}}$ for each $i=0,...,m-1$.
\end{definition}

Clearly, not every Cauchy sequence is $\pi$-local. However, it will be enough for us to work in a set-up where each Cauchy sequence has a $\pi$-local representative in its equivalence class.

\begin{definition}
Suppose that there is a local homeomorphism $\pi:X \to Y$ between two metric spaces $(X,d) $ and $(Y,d')$. Then the metric space $(X,d')$ is said to be {\it $\pi$-locally modelled } on the metric space $(Y,d')$ via $f$ if:
\begin{enumerate} 
    \item Each Cauchy sequence is equivalent to a $\pi$-local Cauchy sequence;
    \item If any two $\pi$-local Cauchy sequences are equivalent, they are strongly equivalent.
    \end{enumerate}
\end{definition}

We will now construct a class of examples. Start with a complex vector space $V$, and take a path connected and semi-locally simply connected open subset $X\subset V$ which satisfies the following conditions:
\begin{enumerate}[label=(LF\arabic*)]
   \item\label{LF1} the closure $\overline{V}\subset X$ is a real manifold with boundary,
   \item\label{LF2} the complement $\overline{V}\setminus V$ is a locally finite union of submanifolds with boundary whose boundary is contained in the boundary of $\overline{V}$.
\end{enumerate}
We will secretly be thinking of the case when $X=V\setminus \Delta$, where $\Delta$ is a locally finite complex hyperplane arrangement, but we do not want to limit ourself to that as our result works in a greater generality. It is easy to show that the following holds: 

\begin{lem}
If a simply connected open $X\subset V$ satisfies (LF1) and (LF2), then it also satisfies
\begin{enumerate}
\item[(LF3)] or each $p \in \overline{X}$ there exists $R>0$ such that for each $0<\epsilon, \epsilon ' < R$ one has:
    \[ \pi_0(B^V_{\epsilon}(p)\cap X)) \cong \pi_0(B^V_{\epsilon'}(p)\cap X) \text{ and } \pi_1(B^V_{\epsilon}(p)\cap X)) \cong \pi_1(B^V_{\epsilon'}(p)\cap X) .\] 
\end{enumerate}
\end{lem}

Let us suppose that $\pi: \widetilde{X} \to X$ is a regular covering space . We need to endow $\widetilde{X}$ with the pullback Riemannian metric: fix a real symmetric non-degenerate bilinear form $g: V \times V \to \mathbb{R}$, which we are going to treat as a constant Riemannian metric, and consider the restriction of $g$ to $X$. We will treat $X$ as a real differentiable manifold, and give $\widetilde{X}$ the differentiable structure compatible with $\pi$. Since the map $\pi$ is a covering space, one might define a pullback Riemannian metric on $\widetilde{X}$ by setting: 
\[ (\pi ^* g)_{x}(v, w) = g_{\pi(x)}(\pi_*(v), \pi_*(w))  \]
for each $x\in \widetilde{X}$ and each $u,v \in T_x \widetilde{X}$. 
Recall that two Riemannian metrics on $\widetilde{X}$ are said to be equivalent if there exist real constants $c$ and $C$ such that 
\[ cg_x(v,v) \leq g'_x(v,v) \leq Cg'_x(v,v) \]
folds for each $x\in X, v\in T_x\widetilde{X}$. To a Riemannian metric we can associate its {\it geodesic distance} as in the previous section, and a Riemannian metric is said to be complete if the induced distance is. Call $\widetilde{d}$ the geodesic distance induced by the pullback Riemannian metric on $\widetilde{X}$, and denote by $\widehat{\widetilde{X}}^g$ the metric completion of $\widetilde{X}$ with respect to this distance. First of all, we need to show that the metric completion is independent on the chosen real inner product on $X$.
\begin{lem}
The topology of the space $\widehat{\widetilde{X}}^g$ does not depend on the metric $g$.
\end{lem}

\begin{proof}
Since any two norms are equivalent on a vector space, one can easily show that any two inner products are equivalent as Riemannian metrics as well. It is straightforward to show that this implies that their respective pullback metrics are equivalent Riemannian metrics, which induce strongly equivalent distances. Any two strongly equivalent distances have the same Cauchy sequences, hence the two metric completions are isomorphic. 

\end{proof}
We will write $\widehat{\widetilde{X}}$, and omit the dependence on the inner product. Our purpose is to show the following result. 

\begin{thm}\label{nice} The metric space $(\widetilde{X}, \widetilde{d})$ is $\pi$-locally modelled on $(X,d)$. 
\end{thm}
First of all, notice that by the Galois correspondence since $\widetilde{X} \to X$ is a covering space, then for each point $\tilde{x} \in \widetilde{X}$ there is an injective homomorphism of groups $\pi _*(\pi _1(\widetilde{X}, \tilde{x})) \to \pi _1(X, \pi(\tilde{x}))  $ which makes $G$ into a normal subgroup of the fundamental group of $X$. Moreover, we recall that the corresponding quotient of the fundamental group acts freely and transitively on the fibers.

We now extend the continuous map $\widehat{\pi}$ to a map $\overline{\pi}: \widehat{\widetilde{X}} \to \widehat{X}=\overline{X}\subset V$ between the two metric completions of the source and target spaces, and we want to obtain a description of the fiber above any point in the boundary of $X$. First of all,
we notice that we can extend the monodromy action to a group action on those fibers: the group $G$ acts by isometries, because the pull-back Riemannian metric is $G$-invariant hence its induced geodesic distance must be, as well, therefore any Cauchy sequence gets transformed by $G$ into another Cauchy sequence whose projection to $X$ tends to the same limiting point. Consider now a point $p\in \Delta  $ and its fiber $\widehat{\pi}^{-1} (p)$. We have the following description.

\begin{lem}\label{discrete}
 The fiber $\widehat{\pi} ^{-1}(p)$ with $p\in \Delta$ is discrete in $\widehat{\widetilde{X}}$.
\end{lem}
\begin{proof}
Suppose by contradiction that there exists a point $x\in$ $\widehat{\pi}^{-1}(p)$ such that for each $\epsilon >0$ the ball $B_{\epsilon}(x)$ in $\widehat{\widetilde{X}}$ meets an infinite number of points of $\widehat{\pi}^{-1}(p)$. Now, this in particular implies that there exists a sequence of points $x_n$ in $\widehat{\pi}^{-1}(p)$ such that $\widetilde{d}(x_n, x) \to 0 $ for $n\to \infty$: one can construct it by simply choosing a point $x_n \in \pi^{-1}(p)$ such that $p \in B_{1/n}(x)\setminus B_{1/n+1}(x)$. For each $n>0$, take any continuous path $\gamma_n: [0,1]\to \widehat{\widetilde{X}}$ such that $\gamma _n(0 )=x$ and $\gamma _n (1)=x_n$ which lives inside $B_{1/n}(x)$. The projection $\pi(\gamma_n)$ is a loop around $p$ in $B^V_{1/n}(p)$  which however is not contractible in $B^V_{1/n}(p)\cap X$ because the original path $\gamma _n$ is the lift of $\pi(\gamma _n)$ and its class only depends on the homotopy class of $\pi(\gamma _n)$. Hence $\pi (\gamma _n)$ must define an element of $\pi _1(B^V_{1/n}(p)\cap X)$ which cannot be in $\pi _1(B^V  _{1/n+1}(p)\cap X)$ since $x_n \not \in B_{1/n+1}(x)$. This show that for each $n>0$ one must have $$\pi_1(B^V_{\epsilon}(p)\cap X)) \not \cong \pi_1(B^V_{\epsilon'}(p)\cap X),$$ hence contradicting \ref{LF1}.
\end{proof}

We recall that every smooth proper effective action on a smooth connected manifold admits a fundamental domain. Indeed, for any $G$-invariant metric $d$ on $M$ and for each point $x\in M$ with trivial stabilizer, one can define the following {\it Dirichlet domain}:
\[ D_x = \{ y\in M \text{ s.t. } d(x,y)<d(gx,y) \text{ for each }g\in G\}. \]
It is standard to show that this is indeed a fundamental domain. In our present case, we let $M=\widetilde{X}$ and $G=\Aut{\pi}$ as above. Then the $G$ action is that of a regular covering space, hence it is smooth, proper and effective. Also, every point has a trivial stabilizer, and our pullback Riemannian metric is $G$-invariant precisely because it is pulled back from the quotient $\widetilde{X}/G = X $. We would like to remind a few properties of the Dirichlet domains obtained with this construction.
\begin{lem}
One has the following for each $x\in \widetilde{X}$:
\begin{enumerate}
    \item The restriction $\pi |_{D_x} : D_x \to \pi(D_x)$ is a homeomorphism.
    \item The image $\pi(D_x)$ is dense in $X$.
\end{enumerate}
\end{lem}
\begin{proof} \
\begin{enumerate}
    \item The covering map $\pi|_{D_x}$ must be an injection: indeed, if $\pi(y)=\pi(y')$, one must have $y'=gy$ for some $g\in G$ because the action of $G$ is transitive on the fibers of $\pi$. But one has that  $D_x \cap gD_x = \emptyset$ unless $g=1$.
    \item By applying $\pi$ one has that:
    \begin{align*}
    \pi (\bigcup _{g\in G }g(\overline{D _x}))&=  \bigcup _{g\in G }\pi(g(\overline{D _x})) \\
    &= \bigcup _{g\in G }\pi(\overline{D _x}) \\
    &= \pi(\overline{D _x}) =X.
\end{align*}
Moreover, consider a point $y \in \overline{D}\subset \widetilde{X}$. Since $\widetilde{X}$ is a metric space, there is a sequence $y_n $ in $D _x$ such that $y_n \to y$. Hence there is a sequence $\pi(y_n) $ in $\pi (D _x)$ which must tend to $\pi (y)$, since one has that
\[ |\pi(y_n) - \pi (y)| \leq \tilde{d}(y_n, y) \to 0 , \]
which yields $\pi(y)\in\overline{\pi(D _x)}$. Therefore $\pi(\overline{D _x}) \subset \overline{\pi(D _x)}$. But since $\pi(\overline{D _x}) = X$, one has that in fact \[ \pi(\overline{D _x}) = \overline{\pi(D _x)} =X . \]
\end{enumerate}
\end{proof}

\begin{proof}[Proof of theorem \ref{nice}]
We need to show that properties (1) and (2) in Definition (\ref{good}) are satisfied.

\begin{enumerate}
 
    \item Choose a point $x\in \widetilde{X}$ and consider the corresponding tiling given by the Dirichlet fundamental domain $D_x$. Consider a Cauchy sequence $\{ \widetilde{y}_n\}\subset \widetilde{X}$ and its projection $\{y_n\} \subset X$. By the previous Lemma, $\pi(D_x)\subset X$ is a dense subset, hence we can assume that $\{ y_n\}$ lies entirely in $\pi(D_x)$, possibly by moving each of its terms in a small neighborhood. Consider now the limit \[\lim _{n\to \infty} \widetilde{y}_n = \widetilde{y} _{\infty} \in \widehat{\widetilde{X}}\] by Lemma \ref{discrete}, the fiber of $\pi(\widetilde{y} _{\infty})=y_{\infty}$ with respect to $\widehat{\pi}$ is discrete, hence we can find $\epsilon >0$ such that the ball of radius $\epsilon$ meets only one point in the special fiber: \[ \overline{ B_{ \epsilon}(\widetilde{y}_{\infty} }) \cap \widehat{\pi} ^ {-1} ( y_{\infty}) = \{ \widetilde{y}_{\infty}\}.\] The projection $B= \pi (B_{ \epsilon}(\widetilde{y}_{\infty} ))$ will be a neighborhood of $y_{\infty}$ in $\overline{X}$ with non-empty interior, hence it will contain infinite terms of the sequence $\{ y_n\}$: since we are only interested in equivalence classes of Cauchy sequences, we can assume without loss of generality that $\{ y_n\}$: lies entirely in $B$. Choose now $g\in G$ such that $gD_x \cap \overline{B_{\epsilon} (\widetilde{y}_{\infty} )} \neq \emptyset$: it exists because the tiles $g\overline{D_x}$, when $g$ varies in $G$, cover $\widetilde{X} $, and because if an open ball meets the boundary of a closed set, then it must also meet the interior. We lift the projection $\{ y_n\}$ of the sequence we started with to the set $gD_x \cap B _{\epsilon}(\widetilde{y}_{\infty} )$,which maps homeomorphically onto its image: the sequence we obtain is $\pi$-local by construction and its limit must live in the closure of the ball $\overline{B_{\epsilon} (\widetilde{y}_{\infty} )}$ and lie over the point $y_{\infty}$. But the only intersection between these two sets consists of the point $\widetilde{y}_{\infty}$, hence this sequence must be equivalent to the sequence $\{\widetilde{y}_n\}$. 
\item Consider two $\pi$-local Cauchy sequences $\{ \widetilde{y}_n\}$ and $\{ \widetilde{y}'_n\}$ in $\widetilde{X}$, and choose a point $x\in \widetilde{X}$. With the same density argument we applied before, we can assume that both lie entirely in $\bigcup gD_x$. Since they are both $\pi$-local, moreover, there must exist $g_1$ and $g_2$ such that $\{ \widetilde{y}_n\}  \subset g_1D_x$ and $\{ \widetilde{y}'_n\} \subset g_2D_x$. By possibly replacing $x$ with $g_1x$, we can furthermore assume that $g_1=1$ and call $g_2=g$. Since $D_x$ is a fundamental domain, we can assume that there exists a collection $g_0,...,g_n$ such that:
\begin{enumerate}
	\item $g_0 = 1$,
	\item $g_n= g$,
	\item $g_ {i+1}=g_ig_1$ and
	\item $g_i\overline{D}_x \cap g_{i+1}\overline{D}_x\neq \emptyset $ for each $i=0,...,n-1$.
\end{enumerate}
We now consider a point $ z\in \overline{D}_x \cap g_1\overline{D}_x$, and take the corresponding fundamental domain $D_z$. It is easy to show that conditions (a), (b) and (c) are satisfied when we replace $D_x$ with $D_z$, and that moreover $g_iD_z$ intersects both $g_iD_x$ and $g_ig_{i+1}D_x$ for each $i$, and viceversa: indeed, the intersection $g_iD_x\cap g_iD_z$ contains the subset $g_i(D_x\cap D_z)$ which is not empty, and a similar fact is true for the intersection $g_iD_z \cap g_ig_1D_x= g_iD_z\cap g_{i+1}D_z$. We consider $\{g^{-1}\widetilde{y}'_n \} \subset D_x$ and look at the partition of $D_x$ given by $(D_z \cap D_x, D_z^c \cap D_x)$. If $ \{ \widetilde{ y}_n\}$ and $\{g^{-1}\widetilde{y}'_n \} \subset D_x$ lie in two different pieces of the partition, then each pair $ ( \widetilde{ y}_n, g^{-1}\widetilde{y}'_n), \ ( g^{-1}\widetilde{y}'_n, 
 g_1\widetilde{ y}_n), \ ( g_1\widetilde{ y}_n, g_1g^{-1}\widetilde{y}'_n), ...$ is a pair of $\pi$-local Cauchy sequences in the same open set $D_x$, $D_z$, $g_1D_x$, $g_1D_z$ etc. If they lie in the same piece, say $D_x \cap D_z$, then their projections both lie in $\pi(D_x \cap D_z)$ hence they can both be lifted to each $g_i(D_x\cap D_z)$ or $g_iD_x \cap g_{i+1}D_z$ homeomorphically. Hence in any of the two case we get a chain of of $\pi$-local strongly equivalent Cauchy sequences relative to the same open set.    
\end{enumerate}
\end{proof}

\section{The pullback metric}\label{pbm}

We will now want to apply our theorem to the present situation. Namely, we go back to the case where $\widetilde { X} = \Stab (\cDD )$, $\pi : \Stab (\cDD) \to X = $ is the central charge map and the image of $\pi$ satisfies conditions (LF1) and (LF2). For each Cauchy sequence $\sigma = \{ \sigma _n\}$ in $\Stab (\cDD) $ we define the following  full subcategory $\mathcal{K} _{\sigma} \subset \mathcal{D}$:

\begin{equation}\label{category}
\mathcal {K}_{\sigma} = \{ E \in \cDD  \ | \ m _{\sigma _n} (E) \to 0 \} . 
\end{equation}

\begin{lem}
The category $K_{\sigma}$ is a thick triangulated subcategory of $\cDD$.
\end{lem}

\begin{proof}
The fact that $\mathcal{K}_{\sigma} $ is a triangulated subcategory follows from the fullness, the fact that $m_{\sigma _n}(E[1]) = m_{\sigma _n}(E)$ for every $E\in \cDD$ and \cite[Proposition 3.1]{Ike}. Now we need to show that it is closed under taking direct summands. This follows from the fact that the Harder-Narasimhan factors of a direct sum are those of both factors, hence: 
\[ m_{\sigma _n}(E\oplus F)= m_{\sigma _n}(E) + m_{\sigma _n}(F). \]
Since these are positive numbers are positive, the fact that $m_{\sigma _n}(E\oplus F)\to 0$ implies that both $m_{\sigma _n}(E),m_{\sigma _n}(F)\to 0$.
\end{proof}
\begin{prop}\label{equal}
If $\{ \sigma _n \} $ and $\{ \tau _n\}$ are $\pi$-local equivalent Cauchy sequences, then $\mathcal{K}_{\sigma } = \mathcal{K} _{\tau}$.
\end{prop}
\begin{proof}
By the previous point, if $\{\sigma _n \}$ and $\{ \tau _n \}$ are equivalent, they are strongly equivalent. First, let us suppose that $\{\sigma _n \}$ and $\{ \tau _n \}$ are very strongly equivalent with respect to some common open set $U$. The fact that $\pi$ restricts to a homeomorphism on $U$ implies that the distances between their slicings must eventually be less than one: indeed if this were not the case, we could look at their images $\pi(\sigma _n)=Z$ and $\pi(\tau _m) = W$ and observe that for $n>N$ they must satisfy 
    \[ ||Z_n-W_n|| < \sin \pi \varepsilon \]
    again since their central charges must form equivalent Cauchy sequences in $V= $ \\ $\Hom _{\mathbb{Z}}(K(D), \mathbb{C})$. This is precisely the condition of \cite[Theorem 7.1]{Bri}, since all the norms on a vector space are equivalent, hence we would be able to find another stability condition with central charge $W_n$ and a slicing with distance less than one from the slicing of $\sigma_n$, which would then have to be in $U$ because their image is in $\pi(U)$, contradicting the injectivity of $\pi$ restricted to $U$. We now suppose that $E\in \cDD$ is such that $m_{\sigma _n}(E) \to 0$, and we want to show that $m_{\sigma _n}(E) \to 0$. We reason by induction on the length of the Harder-Narasimhan filtration of $E$ with resp. to $\sigma_n$ (which stabilizes from some point on since $\{ \sigma _n\} $ is $\pi$-local):
    \begin{enumerate}
        \item First, let us assume that $E \in \mathcal{P}_n(\phi_n) $. Then for each $n$, $E\in \mathcal{Q}_n(\phi _n - \varepsilon _n, \phi _n + \varepsilon _n)$, with $\varepsilon _n <1$, since the distance between the two slicings must be less than one. We know that $|Z_n(E)| =m_{\sigma _n}(E) \to 0$, hence also $|W_n(E)| \to 0$ but $|W_n(E)| \geq (\cos \theta) m_{\tau}(E)$ since $\theta < \pi$, so $m_{\tau}(E)\to 0$.
        \item Now let us look at the last triangle of the Harder-Narasimhan filtration of $E$ with respect to $\sigma _n$:
        \[ E_{m-1}\to E \to F_m \to .\]
        By \cite{Ike}, one has again that $m_{\tau_n}(E) \leq m_{\tau_n}(E_{m-1})+ m_{\tau_n}(F_m)$. The first addendum tends to zero by point (1), while the second tends to zero by the inductive hypotesis since its Harder-Narasimhan filtration has length one less the Harder-Narasimhan filtration of $E$, hence $m_{\tau_n}(E)\to 0.$
\end{enumerate}
\end{proof}
Moreover, we will need the following:

\begin{definition}
We say that a Cauchy sequence of stability conditions $\{\sigma _n = (Z_n, \mathcal{P}_n)\}$ satisfies the {\it limiting support property} if:
\begin{enumerate}
    \item For each $n\in \mathbb{N}$, there exists a $C_n>0$ such that for every $E\in \mathcal{P}_n(\phi)$, $\lim _{n \to \infty }Z_n(E) \neq 0$, one has
    \begin{equation}\label{limsup}
    |Z_n(E)|> C_n ||E||
    \end{equation}
    where $|| \cdot ||$ is any fixed norm on the vector space $\Hom _{\mathbb{Z} }(\Gamma , \mathbb{C})$.
    We will take $C_n$ to be the infimum of the constants which satisfy (\ref{limsup}).
    \item one has, moreover, that 
    \[ \lim _{n\to \infty} C_n = C>0 .\]
\end{enumerate}

\end{definition}

\begin{lem}
The limiting support property is well-defined on equivalence classes of Cauchy sequence of stability conditions.
\end{lem}

\begin{proof}
First, let us consider two Cauchy sequences $\{ \sigma_n=(Z_n, \mathcal{P}_n)\}$ and $\{ \tau _n = (W_n, \mathcal{Q}_n) \}$ in the same equivalence class. Suppose that $\sigma _n$ satisfies the limiting support property with constants $C_n \to C>0$, and $\tau _n $ satisfies the support property with constant $C_n'\to C'>0$. We want to show that:
\[ C - C' = \lim _{n\to \infty} |C_n - C'_n| = 0 . \]
Let us write $|x-y| = d(x,y)$. The fact that $\sigma _n$ and $\tau _n$ are Cauchy sequences implies that also their projections on $\Hom _{\Z} (K(\cDD), \C)$ are, and that moreover they are also in the same equivalence class. Hence $d^*(Z_n, W_n)$ tends to zero with $n$ increasing. This implies that for each $E\in \cDD$ the distance $d(Z_n(E),W_n(E))$ tends to zero as two linear operators coincide if and only if they agree punctually. By a form of the triangle inequality, therefore, one has that \[\left|\frac{|Z(E)|}{||E||}-\frac{|W(E)|}{||E||}\right| \to 0 \] tends to zero for each $E\in \cDD$. Now, since by definition \[ C_n := \inf_{E\in \cDD} \frac{|Z_n(E)|}{||E||} \text{ and } C'_n := \inf_{E\in \cDD} \frac{|W_n(E)|}{||E||}, \] there must be sequences \[ \frac{|Z_n(E_i)|}{||E_i||} \to C_n \text{ and } \frac{|W_n(E_j)|}{||E_j||} \to C'_n \] for each $n$. Hence 
\[ d\left(\frac{|Z_n(E_i)|}{||E_i||}, \frac{|W_n(E_j)|}{||E_j||}\right) \to d\left(\frac{|Z_n(E_i)|}{||E_i||}, C'_n \right)  \text{ and } d\left(\frac{|Z_n(E_i)|}{||E_i||}, C'_n \right) \to d(C_n, C'n) \]
since for any two sequences $a_n\to a$ and $b_n\to b$ one has $d(a_n, b_n) \to d(a,b)$. By uniqueness of limits, this gives $d(C_n, C'_n) \to 0$ and, similarly one has \[ d(C_n, C'_n) \to d(C,C') \to 0.\]
\end{proof}

\section{The limiting slicing}\label{tls}

In order to define the slicing $\overline{\mathcal{P}}_{\infty}$ on the quotient triangulated category $\cDD / \mathcal{K}_{\sigma}$, we need an intermediate step. We will first define a limiting slicing for the Cauchy sequence $\sigma _n$, then we will want to descend it to a slicing on the quotient category. Unfortunately, this is not always possible, or at least non always easy to do. What we will do, instead, is work with hearts of t-structures: as a pair of a central charge and a slicing defines a one-parameter family of t-structures, we will show that at least one of their corresponding hearts can be descended to the quotient category and is moreover compatible with the central charge we have constructed above. Hence we begin with the following:
\begin{definition}\label{limitingslicing}
Let $\{ \sigma _n = (Z_n, \mathcal{P}_n )\}$ be a $\pi$-local Cauchy sequence. Then we define a collection of full, additive subcategories of the triangulated category $\cDD$:
\[ \mathcal{P}_{\infty} = \{ \mathcal{P}_{\infty} (\phi)\}_{\phi \in \mathbb{R}} \]
where
\begin{equation}
\mathcal{P}_{\infty}(\phi) = \{ E\in \cDD \ | \ \phi^{\pm}_{\mathcal{P}_n}(E) \to \phi \} .
\end{equation}
\end{definition}

First of all, we would like to show that such collection of subcategories is a slicing. In order to achieve this result, we need the following:
\begin{prop}\label{walls}
Let ${\sigma _n}$ be a $\pi$-local Cauchy sequence satisfying the limiting support property. Then for each $E\in \cDD$ there exists $N_E>0$ such that the Harder-Narasimhan filtration of $E$ in the stability condition $\sigma _n$ stabilizes after $n> N_E$. 
\end{prop}

\begin{proof}
Fix $E\in \cDD$, $E\not \in \mathcal{K}_{\sigma}$. The HN filtration of $E$ only changes at numerical walls, i.e. when the central charges of two of its semistable factors align on the complex plane, hence if we show that only a finite number of classes in the Grothendieck group $K(\cDD)$ can appear as classes of semistable factors of $E$ for elements of the $\pi$-local Cauchy sequence ${\sigma _n}$, we will always be able to select a subsequence such that only a finite number of numerical walls are crossed. We will do so by bounding the mass of $E$: the fact that $E$ does not belong to $\mathcal{K}_{\sigma}$ implies, by the limiting support property, that its mass $m_{\sigma _n} (E)$ cannot tend to zero when $n\to \infty$. Moreover, since ${\sigma _n}$ is a $\pi$-local Cauchy sequence, the mass of $E$ cannot tend to infinity either: on a small open set which maps homeomorphically onto $X$, the slicings $\{ \mathcal{P}_n\} _{n\in _{\mathbb{N}}}$ of a $\pi$-local Cauchy sequence must form themselves a Cauchy sequence in the standard metric introduced by Bridgeland in \cite[Section 6]{Bri} by the same reasoning as that of Proposition \ref{equal}. So one needs to have $m_{\sigma _n}(E) \to M$ for some $M>0$, hence there is $N>0$ such that for $n>N$ one has $m_{\sigma _n}(E) < 2M$. Suppose that $F$ is a semistable factor of $E$ at a point of the sequence $\sigma _n= (Z_n, \mathcal{P}_n)$ with $n>N$. Then by the support property, one must have \[ ||A|| < \frac{|Z_n(A)|}{C_n} \]. Moreover, since $A$ is a semistable factor of $E$ one must have \[ M_{\sigma _n}(A) = |Z_n(A) | \leq m_{\sigma _n }(E) < 2M\]. Hence we have \[ ||A|| < \frac{|Z_n(A)|}{C_n} \leq \frac{m_{\sigma_n}(E)}{C_n} < \frac{2M}{C_n} \to \frac{2M}{C} . \] This implies that the norm of $A$ must be bounded, hence there can only be a finite number of integral classes which can appear as classes of HN factors for $E$.
\end{proof}

Now we are ready to prove the following:
\begin{prop}
The collection of full additive subcategories defined in (\ref{limitingslicing}) is a slicing on the category $\cDD$.
\end{prop}
\begin{proof}
We need to show that the three properties in the definition of a slicing are satisfied.
\begin{enumerate}
    \item First of all, one has that
    \begin{align*}
        \mathcal{P}_{\infty}(\phi + 1) & = \{ E\in \cDD \ | \ \phi^{\pm}_{\mathcal{P}_n}(E) \to \phi + 1\} \\
        &= \{ E[-1]\in \cDD \ | \ \phi^{\pm}_{\mathcal{P}_n}(E[-1]) \to \phi \} .
    \end{align*} \
    Hence $E\in \mathcal{P}_{\infty}(\phi + 1)$ iff $E[-1]\in \mathcal{P}_{\infty}(\phi)$, and that is only true if $E\in \mathcal{P}_{\infty}(\phi)[1]$.
    \item As far as the Hom-vanishing, we will argue in the following way. We will need to show that:
    \[ \Hom (\mathcal{P}_{\infty}(\phi), \mathcal{P}_{\infty}(\psi) )= 0\]
    whenever $\phi > \psi$. Take $E\in \mathcal{P}_{\infty}(\phi)$ and $F\in \mathcal{P}_{\infty}(\psi)$ and take any $\varepsilon < \frac{\phi - \psi}{2}$. By definition of limit, there must be $N\in \mathbb{N}$ such that for each $n> N$, we get $\phi^{\pm}_{\mathcal{P}_n}(E) \in (\phi-\varepsilon, \phi + \varepsilon)$, respectively $\phi^{\pm}_{\mathcal{P}_n}(F) \in (\psi-\varepsilon, \psi + \varepsilon)$. Hence, $\Hom(E,F)=0$ since $\mathcal{P}_n$ is a slicing for each n.
    \item To construct Harder-Narasimhan filtrations, it is enough to use Proposition (\ref{walls}). Indeed, take $E\in \cDD$. Then the limiting support property implies that any ball centered at $[E]$ in the Grothendieck group $K(\cDD)$ intersects a finite number of walls, hence the Harder-Narasimhan filtration of $E$ with respect to $\sigma _n$ stabilizes for $n$ big. So we can just define  
    \[ \HN ^i_{\infty}(E) =  \HN ^i_n(E) \]
    for $n>>0$.
\end{enumerate}
\end{proof}

\begin{rmk}
Notice that the slicing $\mathcal{P}_{\infty}$ is {\it not} well defined under equivalence of Cauchy sequences: consider the following example. Take $\mathcal{D}=\mathcal{D}_3(A_2)$. This category has a standard heart $\mathcal{A}=\text{Rep}(A_2)$ with two simple generators, $S_2$ and $S_1$, corresponding to the two representations $0\rightarrow \mathbb{C}$ and $\mathbb{C} \rightarrow 0$ respectively, and an irreducible representation $0 \to S_2 \to E \to S_1 \to 0$ corresponding to the representation $\mathbb{C} \overset{1}{\rightarrow}\mathbb{C}$. We will denote a central charge $Z\in \mathbb{C}^2$ as a pair $Z=(Z(S_1), Z(S_2))$. Consider the following two Cauchy sequences: $\sigma _n = (Z_n,\mathcal{P}_n)$ where $Z_n= ( \frac{1}{n}z, -1)$ and the slicing is the only compatible one, and $\sigma _n'= (Z_n', \mathcal{P}_n')$ where $Z_n'=(\frac{1}{n}z', -1)$ with the corresponding slicing and $z, z' \in \mathbb{H}$. Then $\mathcal{P}_{\infty}(1) = \mathcal{P}_{\infty}'(1) = \langle S_2\rangle$, but $\langle S_1 \rangle = \mathcal{P}_{\infty} ( \frac{1}{\pi} \arg z) $ and $\langle S_1 \rangle = \mathcal{P}'_{\infty} ( \frac{1}{\pi} \arg z') $, so $\mathcal{P} \neq \mathcal{P}'$ if $z\neq z'$.

\end{rmk}

The reason why the limiting slicing is not well defined is that it ``remembers" the way objects in $\mathcal{K}$ are sent to zero via the central charge. However, we will see that once we go to the quotient category, this ambiguity will disappear.

Consider now $\mathcal{A} = \mathcal{P}_{\infty}((0,1])$ which, in virtue of the proposition, is the heart of a t-structure. We would like to descend the heart $\mathcal{A}$ to the quotient category $\cDD / \mathcal{K_{\sigma}}$. 

\begin{prop}\label{heart}
The subcategory $\mathcal{A} \cap \mathcal{K} \subset \mathcal{A}$ is a Serre subcategory, hence the t-structure with heart $\mathcal{A} $ descends to a t-structure in $\mathcal{D}/\mathcal{K}$.

\end{prop}

  \begin{proof}
    Consider a short exact sequence $0\to E\to F\to G\to 0$ in $\mathcal{A}$. 
    \begin{enumerate}
        \item First, let us assume that $E$ and $F$ are in $\mathcal{K}$. Then by definition \[m_{\sigma_{n}}(E), m_{\sigma_{n}}(G) \to 0\] if $n\to \infty$. But since $\mathcal{A}$ is the heart of a t-structure in $\cDD$ and therefore the short exact sequence above gives rise to a distinguished triangle, \cite[Proposition 3.3]{Ike} gives $m_{\sigma_{n}}(F)\to 0$ if $n\to \infty$ hence $F\in \mathcal{K} $.
        \item Now, let us assume that $F\in \mathcal{K}$. This implies, in particular, that $Z_{\infty}(F)=0$ since $|Z_{\infty}(F)| \leq \lim _{n\to \infty} m_{\sigma_{n}} (F)=0$. This in turn implies that $|Z_{ \infty}(E)|= -|Z_{\infty}(G)|$ since $Z_{\infty}$ is additive on short exact sequences in $\mathcal{A}$. This would only be possible if some of the Harder-Narasimhan factors of $F$ and/or $G$ have phases tending to zero. But this cannot happen, since the fact that $\mathcal{A}$ is abelian hence closed by taking subobjects and quotients would imply that such Harder-Narasimhan factor, say $H^i$, is in $\mathcal{P}_{\infty}(0)$
        . This implies that $|Z_{\infty} (E)|=|Z_{\infty}(G)|=0$. The reverse triangle inequality gives:
        \begin{align*}
     |Z_{\infty}(E)|= & |\sum Z_{\infty} (\mathrm{HN}^i(E)) | \\
     & \geq \cos \theta \sum |Z_{\infty} (\mathrm{HN}^i(E))| \\
     & = \cos \theta \cdot m_{\sigma _n}(E)
        \end{align*}
        where $\theta$ is the maximum angle between the rays of each $Z_{\infty}(\mathrm{HN}^i(E))$ and that of $Z_{\infty}(E)$, which must be strictly less that $\pi$ because the same reasoning as above gives that all the rays $Z_{\infty}(\mathrm{HN}^i(E))$ have to be contained in the upper half plane $\mathbb{H} \cup \{ 0 \}$. Hence $m_{\sigma _n}(E) \to 0$, which gives $E\in \mathcal{K}$, and the same is true for $G$.
        
    \end{enumerate}
 The fact that this implies that the t-structure with heart $\mathcal{A}$ descends to a t-structure on $\mathcal{D}/\mathcal{K}$ follows from \cite[Proposition 2.20]{AGH}.
    \end{proof}
    Moreover, we have the following:
    \begin{lem}
For any two equivalent $\pi$-local Cauchy sequences the induced hearts on the quotient category $\cDD / \mathcal{K}$ coincide.
\end{lem}
\begin{proof}
Take two equivalent $\pi$-local Cauchy sequences $\sigma _n$ and $\sigma_n'$. The limit of this sequence will be $\sigma _{\infty}= \sigma _{\infty}'=(Z_{\infty}, \mathcal{P}_{\infty})$. Consider $E\in \mathcal{P}_{\infty}(\phi)$ for some $\phi \in \mathbb{R}$. Now, either $Z_{\infty}(E)\neq 0$ or $Z_{\infty}(E)= 0$. If $Z_{\infty}(E)\neq 0$, then for compatibility reasons one must have $\phi = \frac{1}{\pi}\text{arg}Z_{\infty}(E)$. If $Z_{\infty}(E)= 0$, then it must be $E\in \mathcal{K}_{\sigma}=\mathcal{K}_{\sigma'}$, since $E\in \mathcal{P}_{\infty}(\phi)$ implies $E\in \mathcal{P}_n(\phi - \epsilon, \phi + \epsilon)$ for $n>>0$, hence its limiting central charge cannot vanish without its whole mass vanishing. Hence the only ambiguity in the limiting slicing happens precisely for objects in $\mathcal{K}$, which implies that $\mathcal{P}_{\infty}((0,1])/ \mathcal{P}_{\infty}((0,1])\cap \mathcal{K}$ does not depend on the $\pi$-local representative in the equivalence class.
\end{proof}

\section{The locus of deformable degenerate stability conditions}\label{lddsc}
    We can now define our map $j$ of Theorem \ref{Main}: take an equivalence class of Cauchy sequences $[\sigma]$. First of all, the content of Theorem \ref{nice} is that we can always choose a $\pi$-local representative. Suppose our class is $[\sigma ] \in \widehat{\Stab}(\cDD)$ and that $\sigma$ is a $\pi$-local representative. We define:
    
\[ j(\underline{\sigma} )= (\mathcal{K}_{\sigma}, (\overline{Z}_{\infty}, \overline{\mathcal{A}})) \]

where:
\begin{enumerate}
\item The category $\mathcal{K}_{\sigma} $ is as in \ref{category}. Notice that Proposition \ref{equal} implies that $\mathcal{K}_{\sigma}$ is well defined because two $\pi$-local equivalent Cauchy sequences give the same subcategories.
    \item The group homomorphism $\overline{Z}_{\infty}: K(\cDD / \mathcal{K}) \to \mathbb{C}$ is induced by $Z_{\infty} = \lim _{n\to \infty} Z_n$ using the short exact sequence 
    \[ K(\mathcal{K}) \to K(\mathcal{D}) \to K(\mathcal{D} / \mathcal{K} ) \to 0 \]
    considering that $Z_{\infty}|_{K(\mathcal{K})} = 0$.
    \item The heart $\overline{\mathcal{A}}$ in $\cDD / \mathcal{K}$ is the heart of the induced t-structure of Proposition \ref{heart}.

\end{enumerate}

\

First, notice that we have the following:
\begin{thm}
The pair $(\overline{Z}_{\infty}, \overline{\mathcal{A}})$ defines a stability condition on $\mathcal{D}/\mathcal{K}$.
\end{thm}

\begin{proof}
We need to show that:
\begin{enumerate}
    \item If $E\in \overline{\mathcal{A}}$, then $\overline{Z_{\infty}}(E)\in \mathbb{H}$;
    \item There exists a finite HN filtration for each object;
    \item The support property is satisfies.
\end{enumerate}
For (1), take $[E]\in \overline{\mathcal{A}}$. First of all, notice that $Z_{\infty}$ is well defined on $\mathcal{D}/\mathcal{K}$: any two liftings $E$ and $E'$ of $[E]$ differ by elements of $\mathcal{K}$ in the sense that there is a morphism $f: E\to E'$ such that $\ker (f) $ and $\mathrm{coker}(f)$ both lie in $\mathcal{K}$. By looking again at the two short exact sequences $0 \to \ker (f) \to E \to \mathrm{im}(f) \to 0 $ and $0 \to \mathrm{im}(f) \to E' \to \mathrm{coker}(f) \to 0$ one concludes that 

\begin{align*}
    Z_{\infty} (E) = & Z_{\infty} (\ker (f)) + Z_{\infty}(\mathrm{im}(f)) \\
    & = Z_{\infty}(\mathrm{im}(f)) \\
    & = Z_{\infty}(E') - Z_{\infty}(\mathrm{coker}(f))\\
    & = Z_{\infty} (E').\\
    \end{align*}
Choose a lifting $E$ of $[E]$. By the same reasoning as Proposition \ref{heart}, one shows that $Z_{\infty}(E)\in \mathbb{H}\cup \{ 0\}$. Therefore we only need to exclude that $Z_{\infty}(E)=0$: indeed, if that is the case, we have already shown that this must imply $m_{\sigma _{\infty}}(E) =0$, hence $E\in \mathcal{K}$ and so $[E]=[0]$.\\
For (2) it is enough to take the image of the Harder-Narasimhan filtration of an object via the quotient functor $p: \mathcal{D} \to \mathcal{D}/\mathcal{K}$ since the heart $\mathcal{A}$ has the Harder-Narasimhan property. \\
For (3), this discends directly from the fact that we have assumed our sequence to have the limiting support property: indeed, we know that
\[ \text{inf}\left\{ \frac{|Z _{\infty}(E)|}{||E||} \ \Bigg| \ E \in \mathcal{A}\text{ } Z_{\infty}-\text{semistable}, E \not \in \mathcal{K}\right\} = C>0 \]
for a suitable constant $C$ and a fixed norm $|| \cdot ||$ on $K(\mathcal{D}) \otimes \mathbb{C}$. Now, we fix a special norm on $K(\mathcal{D}/\mathcal{K}) \otimes \mathbb{C}$ which is a quotient of $K(\mathcal{D} )\otimes \mathbb{C}$ as we have seen before: we let
\[ ||[E]||_q = \text{inf} \{ || E + F || \ | \ F\in K \} .\]

It is a classical exercise to show that this is a norm. This, in particular, implies that $||[E]||_q \leq ||E||$. Note that an object is $\overline{Z}$-semistable in $\overline{\mathcal{A}}$ if and only if its Harder-Narasimhan factors $\text{HN}^1(E), \ldots, \text{HN}^{k-1}(E)$ but either the first or the last one are in $\mathcal{K}$. Hence we have:

\begin{align*}
    \text{inf}\left\{ \frac{|\overline{Z}_{\infty}([E]) |}{||[E]||_q}  \ \Bigg| \ E \in \overline{\mathcal{A}}\text{ }\overline{ Z }_{\infty}-\text{ss},\right\}& = \text{inf}\left\{ \frac{|Z_{\infty}(E) |}{||[E]||_q} \ \Bigg| \ E \in \overline{\mathcal{A}}\text{ }\overline{ Z }_{\infty}-\text{ss}\right\}\\ &=\text{inf} \left\{ \frac{|Z_{\infty}(\text{HN}^i(E)) + \sum _{j\neq i} Z _{\infty}(\text{HN}^j(E))|}{||[E]||_q} \ \Bigg| \ E \in \overline{\mathcal{A}}\text{ }\overline{ Z }_{\infty}-\text{ss}\right\} \\ & =\text{inf} \left\{ \frac{|Z_{\infty}(\text{HN}^i(E)))|}{||[\text{HN}^i(E)]||_q}\ \Bigg| \ \text{HN}^i(E) \in \mathcal{A}\text{ } Z _{\infty}-\text{ss}\right\} \\ & \geq \text{inf} \left\{ \frac{|Z_{\infty}(E') |}{||E'||} \ \Bigg| \ E' \in \mathcal{A} \ Z_{\infty}-\text{ss}\right\} =C>0 .
\end{align*}
\end{proof}

To end the proof of Theorem \ref{Main}, moreover, we have the following:

\begin{prop}
The map $j$ is injective.
\end{prop}

\begin{proof}
We need to show that, if we take any two $\pi$-local Cauchy sequences $\underline{\sigma}=\{\sigma _n \} = \{(Z_n, \mathcal{P}_n)\}$ and $\underline{\tau}=\{\tau _n\}= \{ (W_n, \mathcal{Q}_n)\}$ such that $j(\underline{\sigma})=j(\underline{\tau})$, then it must be $\sigma \sim \tau$, i.e. $\tilde{d}(\sigma_n, \tau_n) \to 0$ as $n \to \infty$. Suppose that $j(\underline{\sigma})=(\mathcal{K}_{\sigma}, \underline{Z}, \underline{\mathcal{P}})$ and $j(\underline{\tau})=(\mathcal{K}_{\tau}, \underline{W}, \underline{\mathcal{Q}})$ with $\mathcal{K}_{\sigma} = \mathcal{K}_{\tau}$ and $\underline{Z}=\underline{W}$ but $\underline{\sigma}\not \sim \underline{\tau}$. We want to show that $\underline{\mathcal{P} }\neq \underline{\mathcal{Q}}$. First of all, we claim that if $d(\mathcal{P}_{\infty}, \mathcal{Q}_{\infty}) \geq 1$, where the distance is the usual distance on the set of slicings
\[ d(\mathcal{P}, \mathcal{Q}) = \sup _{E\in \mathcal{D}} \{ |\phi^-_{\mathcal{P}}(E) -  \phi^-_{\mathcal{Q}}(E) |, | \phi^+_{\mathcal{P}}(E) - \phi^+_{\mathcal{Q}}(E) | \} \]
as defined in  \cite{Bri}, and there exists $E\not \in \mathcal{K}= \mathcal{K}_{\sigma}= \mathcal{K}_{\tau}$ such that 
\begin{equation}\label{sup}
 \sup \{ |\phi^-_{\mathcal{P}}(E) -  \phi^-_{\mathcal{Q}}(E) |, | \phi^+_{\mathcal{P}}(E) - \phi^+_{\mathcal{Q}}(E) | \} \geq 1 
 \end{equation}
then $\overline{\mathcal{P}} \not = \overline{\mathcal{Q}} $. Indeed, if this is the case, then up to shifting $E$ we can assume that $E\in \mathcal{P}_{\infty}((0,1]) $ but $E\not \in \mathcal{Q}_{\infty}((0,1])$, which if $E\not \in \mathcal{K}$ implies that $E\in \overline{\mathcal{A}} = \mathcal{P}_{\infty} ((0,1]) /  \mathcal{P}_{\infty} ((0,1]) \cap \mathcal{K} $ but $E\not \in  \mathcal{Q}_{\infty} ((0,1]) /  \mathcal{Q}_{\infty} ((0,1]) \cap \mathcal{K} $. Therefore, since the induced central charges coincide, the induced slicings must be different. Hence we need to prove that there is $E\in \mathcal{D}, E\not \in \mathcal{K}$ such that \ref{sup} is achieved. The fact that both the thick subcategory $\mathcal{K}$ and the induced central charges on the quotient triangulated category coincide, the fact that $\overline{\sigma} \not \sim \overline{\tau}$ implies that the respective limiting points in the metric completion must bot lie in the same fiber. By proposition \ref{discrete} the fiber is discrete, hence we can isolate the two limiting points $\sigma _{\infty}$ and $\tau _{\infty}$ with balls $B$ and $B'$ in $\widehat{\Stab(\mathcal{D})}$. With the same method as that of Theorem \ref{nice}, we can produce a sequence $\underline{\tau '}$ which is equivalent to $\underline{\sigma}$ and such that $\underline{\tau '} = \{(W_n, \mathcal{Q}'_n)\}$, possibly up to tweaking the sequence $\underline{\tau}$ slightly. Now by \cite[Proposition]{Bri}, for each $n>0$ one must have $d(\mathcal{Q}_n,\mathcal{Q}_n')\geq 1$, otherwise the sequences $\underline{\tau}$ and $\underline{\tau'}$ would coincide pointwise, thus violating the non-equivalence of $\underline{\sigma}$ and $\underline{\tau}$. Moreover, it must be $\mathcal{Q}'=\mathcal{Q}[k]$, for some $k\in \mathbb{ Z}, k\neq 0$ which implies that the inequality \ref{sup} is achieved for each $E\in \mathcal{D}$. Hence their limits must have distance $d(\mathcal{Q}_{\infty}, \mathcal{Q}_{\infty})\geq 1$. Moreover, since $\sigma$ and $\tau'$ are now $\pi$-local equivalent Cauchy sequences, one must have that $d(\mathcal{P}_n, \mathcal{Q}'_n )\to 0$. This implies that 
\[ d(\mathcal{P}_n, \mathcal{Q}_n) + d(\mathcal{P}_n, \mathcal{Q}_n') \geq d(\mathcal{Q}_n, \mathcal{Q}_n') \to d(\mathcal{P}_{\infty}, \mathcal{Q}_{\infty}) ) \geq d(\mathcal{Q}_{\infty}, \mathcal{Q}'_{\infty}) \geq 1\]
and the inequality \ref{sup} must be achieved for some $E\not \in \mathcal{K}$.
\end{proof}

\begin{exam}
We consider the quiver $Q=A_1=\{ \bullet \}$, and take the category $\mathcal{D}_N = \mathcal{D}(\Gamma _N(Q))$. Such category is generated by a single spherical object $S$ corresponding to the constant path at the only vertex of $Q$. Clearly, one has that $K(\mathcal{D})\cong \mathbb{Z}$, hence the space of central charges is $\mathbb{C}^* \subset \Hom _{\mathbb{Z}}(\mathbb{Z}, \mathbb{C})\cong \mathbb{C}$, and that the restriction of the forgetful map onto its image \[ \pi : \Stab (\cDD ) \cong \widetilde{\mathbb{C}^*} \to \mathbb{C}^* \] is the universal covering map. Now we need to look at two spaces:
\begin{enumerate}
    \item The metric completion $\overline{\widetilde{\mathbb{C}^*}}$ with respect to the pullback of the restriction of the Euclidean metric from $\mathbb{C}^*$. Such space is, set-theoretically, \[ \widetilde{\mathbb{C}^*} \cup * = \mathbb{C} \cup * . \] Indeed, take a sequence $\{ \sigma _n \} \in \widetilde{\mathbb{C}^*}$. Such sequence is Cauchy with respect to the distance $d$ induced by the pullback metric if 
    \[ d(\sigma _n, \sigma _m) = \inf \left\{ \int _0 ^1 ||\gamma'(t)|| \ | \ \gamma(0)=\sigma_n, \gamma(1)=\sigma _m \right\} \to 0\] as $m,n\to \infty$. Since \[ ||\gamma'(t)|| = \langle \pi_* \gamma'(t), \pi_* \gamma'(t) \rangle _{Eu} ^{ 1/2 }, \] such infimum would be achieved by a path projecting onto a geodesic in $\mathbb{C}^*$, i.e. by the path going directly from $\pi(\sigma _n)$ to $\pi(\sigma _m)$ if $\sigma_n$ and $\sigma_m$ are on the same sheet of the universal cover, or by the path going through the (missing) preimage of zero which projects onto the union of the two segments connecting $\pi(\sigma _n)$ and $\pi(\sigma _n)$ to zero. Hence one sees that if one takes two points $\sigma _n = \log z$ and $\sigma _m = \log z + 2k\pi$ their distance is 
    \[ d( \log z, \log z + 2k\pi ) = |z| , \] which implies that any Cauchy sequence is equivalent to a Cauchy sequence consisting entirely on points lying on one sheet of the universal cover, which is isometric to $\mathbb{C}^*$ via the restriction of the projection. Hence the boundary $\overline{\Stab}(\cDD) \setminus \Stab(\cDD)$ consists only of a point.
    \item The set $\Stab _{sd}(\cDD)$. Since the object $S$ is simple, there are no proper thick subcategories of $\cDD$. Hence we have two cases: either $\mathcal{K}=(0)$, or $\mathcal{K}=\cDD$. In the first case the Verdier localization gives $\cDD / (0) = \cDD$, while in the second case one has $\cDD / \cDD = (0)$. Hence the set $\Stab _{sd}(\cDD)$ is again in bijection with the set \[\Stab (\cDD ) \cup * \cong \widetilde{\mathbb{C}^*} \cup * .\] 
    \end{enumerate}
    By choosing a branch of the complex logarithm and by identifying a central charge $Z$ with $\log z= \log Z(S) \in \mathbb{C}$, we can define our map
    \[ j: \overline{\Stab(\cDD) } \longrightarrow \Stab(\cDD)_{sd}, \] which in this case is bijective, as:
    \[ j(([\log z_n], \mathcal{A}) )=  \left\{\begin{array}{lr}
        (\mathcal{K}=(0), (\lim _{n\to \infty}\log z_n, \mathcal{A})) & \text{for } \lim _{n\to \infty}\log z_n \neq 0\\
        (\mathcal{K}=\mathcal{\cDD}, (0, 0)) & \text{for } \lim _{n\to \infty}\log z_n = 0\\
        \end{array}\right\} .\]
\end{exam}

\end{document}